\documentclass[10pt]{article}%

\usepackage{rotating} 

\usepackage{mathrsfs}
\DeclareMathAlphabet{\mathpzc}{OT1}{pzc}{m}{it}

\newcommand{\starA}{\mathscr{A}_k}

\usepackage{color}
\usepackage{amsmath}
\usepackage{graphicx}
\usepackage{amsfonts}
\usepackage{amssymb}%
\usepackage{latexsym}
\usepackage{psfrag}
\usepackage{subfigure}
\usepackage{accents}
\usepackage{bm}

\newtheorem{theorem}{Theorem}[section]
\newtheorem{corollary}[theorem]{Corollary}
\newtheorem{definition}[theorem]{Definition}
\newenvironment{proof}[1][Proof]{\noindent \emph{#1.} }
{\hfill \ \rule{0.5em}{0.5em}}
\newtheorem{lemma}[theorem]{Lemma}
\newtheorem{proposition}[theorem]{Proposition}

\newtheorem{assumption}[theorem]{Assumption}
\newtheorem{example}[theorem]{Example}
\newtheorem{remark}[theorem]{Remark}

\numberwithin{equation}{section}

\usepackage[a4paper]{geometry}
\newcommand{\noi}{\noindent}

\newcommand{\R}{\mathbb{R}}

\newcommand{\cV}{{\cal V}}

\newcommand{\cO}{{\cal O}}

\newcommand{\bx}{\mathbf{x}}
\newcommand{\br}{\mathbf{r}}
\newcommand{\bn}{\mathbf{n}}

\newcommand{\ba}{\mathbf{a}}
\newcommand{\bA}{\mathbf{A}}
\newcommand{\bv}{\mathbf{v}}
\newcommand{\bw}{\mathbf{w}}

\newcommand{\by}{\mathbf{y}}

\newcommand{\bff}{\mathbf{f}}

\newcommand{\bze}{\mathbf{0}}

\newcommand{\sspan}{\mathop{{\rm span}}}

\newcommand{\re}{{\rm e}}
\newcommand{\ri}{{\rm i}}
\newcommand{\rd}{{\rm d}}

\newcommand{\beq}{\begin{equation}}
\newcommand{\eeq}{\end{equation}}
\newcommand{\beqs}{\begin{equation*}}
\newcommand{\eeqs}{\end{equation*}}
\newcommand{\bit}{\begin{itemize}}
\newcommand{\eit}{\end{itemize}}
\newcommand{\ben}{\begin{enumerate}}
\newcommand{\een}{\end{enumerate}}
\newcommand{\bal}{\begin{align}}
\newcommand{\eal}{\end{align}}
\newcommand{\bals}{\begin{align*}}
\newcommand{\eals}{\end{align*}}
\newcommand{\bse}{\begin{subequations}}
\newcommand{\ese}{\end{subequations}}
\newcommand{\bpr}{\begin{proposition}}
\newcommand{\epr}{\end{proposition}}
\newcommand{\bre}{\begin{remark}}
\newcommand{\ere}{\end{remark}}
\newcommand{\bpf}{\begin{proof}}
\newcommand{\epf}{\end{proof}}
\newcommand{\ble}{\begin{lemma}}
\newcommand{\ele}{\end{lemma}}
\newcommand{\bco}{\begin{corollary}}
\newcommand{\eco}{\end{corollary}}
\newcommand{\bex}{\begin{example}}
\newcommand{\eex}{\end{example}}
\newcommand{\bth}{\begin{theorem}}
\newcommand{\enth}{\end{theorem}}

\newcommand{\Rea}{\mathbb{R}}
\newcommand{\Com}{\mathbb{C}}

\newcommand{\essinf}{\mathop{{\rm ess} \inf}}

\newcommand{\Oi}{{\Omega}}

\newcommand{\Oe}{{\Omega_+}}

\newcommand{\eps}{\varepsilon}

\newcommand{\pdiff}[2]{\frac{\partial #1}{\partial #2}}

\newcommand{\dnpu}{\partial_n^+ u}

\newcommand{\gpu}{\gamma^+ u}

\newcommand{\gu}{\nabla u}

\newcommand{\nT}{\nabla_{\bound}}

\newcommand{\half}{\frac{1}{2}}

\newcommand{\bound}{{\partial \Oi}}

\newcommand{\LtG}{{L^2(\bound)}}

\newcommand{\LtGt}{{\LtG\rightarrow \LtG}}
\newcommand{\LtGtH}{{\LtG\rightarrow \HoG}}

\newcommand{\HoG}{H^1(\bound)}
\newcommand{\HohG}{H^1_k(\bound)}

\newcommand{\tendi}{\rightarrow \infty}
\newcommand{\tendo}{\rightarrow 0}

\newcommand{\intotp}{\int^{2\pi}_{0}}

\newcommand{\opA}{A'_{k,\eta}}
\newcommand{\opABW}{A_{k,\eta}}
\newcommand{\opAinv}{(A'_{k,\eta})^{-1}}

\newcommand{\normAinv}{\|\opAinv\|}

\newcommand{\DtN}{P_{\rm{DtN}}}
\newcommand{\NtD}{P_{\rm{NtD}}}

\def\XXint#1#2#3{{\setbox0=\hbox{$#1{#2#3}{\int}$}
     \vcenter{\hbox{$#2#3$}}\kern-.5\wd0}}

\usepackage{hyperref}
\usepackage{tikz}
\definecolor{myblue}{rgb}{0,0,0.6}
\hypersetup{colorlinks=true,
linkcolor=myblue,citecolor=myblue,filecolor=myblue,urlcolor=myblue}
\newcommand*{\N}[1]{\left\|#1\right\|}

\allowdisplaybreaks[4]

\newcommand{\tfa}{\text{ for all }}
\newcommand{\tfor}{\text{ for }}

\newcommand{\tas}{\text{ as }}
\newcommand{\tand}{\text{ and }}

\newcommand{\vertiii}[1]{{\left\vert\kern-0.25ex\left\vert\kern-0.25ex\left\vert #1 
    \right\vert\kern-0.25ex\right\vert\kern-0.25ex\right\vert}}

\newcommand{\Hna}{H_{n}^{(1)}(k )}
\newcommand{\Hnpka}{H_{n}^{(1)'}(k )}

\newcommand{\Gammaext}{\widetilde{\Gamma}}

\usepackage{chngcntr}

\graphicspath{{Figures/}}

\renewcommand{\Re}{\mathop{\rm Re}\nolimits}

\newcommand{\bl}{\begin{flushleft}}
\newcommand{\el}{\end{flushleft}}
\newcommand{\ert}{\end{flushright}}
\newcommand{\bc}{\begin{center}}
\newcommand{\ec}{\end{center}}

\newcommand{\numList}{\begin{enumerate}}
\newcommand{\enumList}{\end{enumerate}}

\newcommand{\Sk}{S_k}
\newcommand{\Dkreg}{C^{2,\alpha}}

\begin{document}

\title{Wavenumber-explicit analysis for the Helmholtz $h$-BEM: error estimates and iteration counts for the Dirichlet problem}
\author{Jeffrey Galkowski\footnotemark[1]\,,\, Eike H.~M\"uller\footnotemark[2]\,\,, Euan A.~Spence\footnotemark[3]}

\date{\today}

\renewcommand{\thefootnote}{\fnsymbol{footnote}}

\footnotetext[1]{Department of Mathematics, Stanford University, Building 380, Stanford, California 94305, USA \tt jeffrey.galkowski@stanford.edu}
\footnotetext[2]{Department of Mathematical Sciences, University of Bath, Bath, BA2 7AY, UK, \tt E.Mueller@bath.ac.uk}
\footnotetext[3]{Department of Mathematical Sciences, University of Bath, Bath, BA2 7AY, UK, \tt E.A.Spence@bath.ac.uk }

\maketitle

\begin{abstract}
We consider solving the exterior Dirichlet problem for the Helmholtz equation with the $h$-version of the boundary element method (BEM) using the standard second-kind combined-field integral equations.
We prove a new, sharp bound on how the number of GMRES iterations must grow with the wavenumber $k$ to have the error in the iterative solution bounded independently of $k$ as $k\tendi$
when the boundary of the obstacle is analytic and has strictly positive curvature.
 To our knowledge, this result is the first-ever sharp bound on how the number of GMRES iterations depends on the wavenumber for an integral equation used to solve a scattering problem.
We also prove new bounds on how $h$ must decrease with $k$ to maintain $k$-independent quasi-optimality of the Galerkin solutions as $k \tendi$
when the obstacle is nontrapping.

\

\noi \textbf{Keywords:} Helmholtz equation,  high frequency, boundary integral equation, boundary element method, GMRES, pollution effect, semiclassical

\

\noi\textbf{AMS Subject Classifications:} 35J05, 35J25, 65N22, 65N38, 65R20

\end{abstract}

\renewcommand{\thefootnote}{\arabic{footnote}}

\section{Introduction}

This paper is concerned with the wavenumber-explicit numerical analysis of boundary integral equations (BIEs) for the Helmholtz equation
\beq\label{eq:Helm}
\Delta u+k^2 u=0, 
\eeq
where $k>0$ is the \emph{wavenumber}, posed in the exterior of a 2- or 3-dimensional bounded obstacle $\Oi$
with Dirichlet boundary conditions on $\Gamma:=\partial \Oi$. 

We consider the standard second-kind combined-field integral equation formulations of this problem: the so-called ``direct" formulation (arising from Green's integral representation)
\beq\label{eq:direct}
\opA v = f_{k,\eta}
\eeq
and the so-called ``indirect" formulation (arising from an ansatz of layer potentials not related to Green's integral representation)
\beq\label{eq:indirect}
\opABW \phi = g_k,
\eeq
where 
\beq\label{eq:scpo}
\opA := \half I + D'_k -\ri\eta S_k, \qquad \opABW := \half I + D_k -\ri\eta S_k,
\eeq
$\eta \in \Rea\setminus \{0\}$ is an arbitrary coupling parameter, $S_k$ is the single-layer operator, $D_k$ is the double-layer operator, and $D'_k$ is the adjoint double-layer operator  \eqref{eq:SD'}, \eqref{eq:D}.

For simplicity of exposition, we focus on the direct equation \eqref{eq:direct}, but the main results also hold for the indirect equation \eqref{eq:indirect} (see Remark \ref{rem:indirect} below).
The contribution to Equation \eqref{eq:direct} from the Dirichlet boundary conditions is contained in the right-hand side $f_{k,\eta}$; our results are independent of the particular form of $f_{k,\eta}$, and so we can simplify the presentation by restricting attention to the particular 
exterior Dirichlet problem corresponding to scattering by a point source or plane wave, i.e.~the \emph{sound-soft scattering problem} (Definition \ref{def:SSSP} below).

We consider solving the equation \eqref{eq:direct} in $\LtG$ using
the Galerkin method; this method seeks an approximation $v_N$ to the solution $v$ from a finite-dimensional approximation space $\cV_N$ (where $N$ is the dimension, i.e.~the total number of degrees of freedom).
In the majority of the paper $\partial \Omega$ is $C^2$, in which case $\cV_N$ will be the space of 
 piecewise polynomials of degree $p$, for some fixed $p\geq 0$, on shape-regular meshes of diameter $h$, with $h$ decreasing to zero;
this is the so-called \emph{$h$--version} of the Galerkin method, and we denote $\cV_N$ and $v_N$ by $\cV_h$ and $v_h$, respectively, and note that $N\sim h^{-(d-1)}$, where $d$ is the dimension. 
To find the Galerkin solution $v_h$, one must solve a linear system of dimension $N$; in practice this is usually done using Krylov-subspace iterative methods such as the generalized minimal residual method (GMRES).

For the numerical analysis of this situation when $k$ is large, there are now, roughly speaking, two main questions:
\ben
\item[Q1.] How must $h$ decrease with $k$ in order to maintain accuracy of the Galerkin solution as $k\tendi$?
\item[Q2.] How does the number of GMRES iterations required to achieve a prescribed accuracy grow with $k$?
\een
The goal of this paper is to prove rigorous results about these two questions, and then compare them with the results of numerical experiments.

We now give short summaries of the main results. These results depend on the choice of the coupling parameter $\eta$; for the results on Q1 we need $|\eta|\sim k$ and for the results on Q2 we need $\eta \sim k$, where we use the notation $a\sim b$ to mean that there exists $C_1, C_2>0$, independent of $h$ and $k$,  such that $C_1 b\leq a \leq C_2 b$. We also use the notation $a\lesssim b$ to mean that there exists $C>0$, independent of $h$ and $k$, such that $a\leq Cb$.

\paragraph{Summary of main results regarding Q1 and their context.} 
Numerical experiments indicate that, in many cases, the condition $hk\lesssim 1$ is sufficient for the
Galerkin method to be quasi-optimal (with the constant of quasi-optimality independent of $k$; i.e., \eqref{eq:qothm} below holds); see \cite[\S5]{GrLoMeSp:15}. This feature can be described by saying that the $h$-BEM does not suffer from the pollution effect (in constrast to the $h$-FEM; see, e.g., \cite{BaSa:00}, \cite[Chapter 4]{Ih:98}). The best existing result in the literature is that $k$-independent quasi-optimality of the Galerkin method applied to the integral equation \eqref{eq:direct} holds when $hk^{(d+1)/2}\lesssim 1$ for 2- and 3-d $C^{2,\alpha}$ obstacles that are star-shaped with respect to a ball \cite[Theorem 1.4]{GrLoMeSp:15}. In this paper we improve this result by showing that the $k$-independent quasioptimality holds for 2-d \emph{nontrapping} obstacles when $hk^{3/2}\lesssim 1$, for 3-d nontrapping obstacles when $hk^{3/2}\log k \lesssim 1$, and for 2- and 3-d smooth (i.e.~$C^\infty$) convex obstacles with strictly positive curvature when $hk^{4/3}\lesssim 1$ (see Theorem \ref{thm:1} below).

The ideas behind the proofs of these results are summarised in Remark \ref{rem:ideaQ1} below, but we highlight here that all the integral-operator bounds used in these arguments are sharp up to a factor of $\log k$. 
Therefore, to lower these thresholds on $h$ for which quasi-optimality is proved, one would need to use different arguments than in the present paper.
We also highlight that recent experiments by 
 Baydoun and Marburg \cite{Ma:16, BaMa:17, Ma:17, BaMa:18} give examples of Helmholtz problems where the $h$-BEM suffers from a pollution effect, and therefore 
determining the sharp threshold on $h$ for $k$-independent quasi-optimality to hold in general
 is an exciting open question.

\paragraph{Summary of main results regarding Q2 and their context.} There has been a large amount of research effort expended on understanding empirically how iteration counts for integral-equation formulations of scattering problems involving the Helmholtz or Maxwell equations depend on $k$; see, e.g, 
\cite{AlBoLe:07, AnDa:07, BoBrLeTu:15, BrElPaTu:09, ViGrGi:14}, 
and the references therein.

To our knowledge, however, there are no sharp $k$-explicit bounds in the literature,
for any integral-equation formulation of a Helmholtz or Maxwell scattering problem,
 on the number of iterations GMRES requires to achieve a prescribed accuracy. 
The main reason, in this current setting of the Helmholtz exterior Dirichlet problem, is that the operator $\opA$
 is non-normal for all obstacles other than the circle and sphere \cite{BeSp:11, BePhSp:13}. Therefore, for sufficiently-accurate discretisations, the Galerkin matrix of $\opA$ is also non-normal, and one cannot use the well-known bounds on GMRES iterations in terms of the condition number (see, e.g., the review in \cite[\S6]{SiSz:07}).

In this paper, we prove that, for 2- and 3-d analytic obstacles with strictly positive curvature, the number of GMRES iterations growing like $k^{1/3}$ is sufficient to have the error in the iterative solution bounded independently of $k$ (see Theorem \ref{thm:GMRES1} below). Numerical experiments in \S\ref{sec:num} show that the numbers of GMRES iterations for the sphere and an ellipsoid grow slightly less than $k^{1/3}$ ($k^{0.29}$ for the sphere and $k^{0.28}$ for an ellipsoid), and thus our bound is effectively sharp.

The ideas behind the proof are summarised in Remark \ref{rem:ideaQ2} below.
The focus of this paper is in proving results for the operator $\opA$, i.e.~the operator in the standard second-kind integral formulation, but we highlight in Remark \ref{rem:scom} below how a bound on the number of GMRES iterations of $k^{1/2}$ when $d=2$ and $k^{1/2}\log k$ when $d=3$ can be obtained for a modification of $\opA$, the so-called \emph{star-combined integral equation} introduced in \cite{SpChGrSm:11}. Moreover, whereas our bound on the number of iterations of $k^{1/3}$ for $\opA$ holds for analytic obstacles with strictly positive curvature, the bounds for the star-combined operator hold for a much wider class of obstacles, namely piecewise-smooth Lipschitz obstacles that are star-shaped with respect to a ball.

\paragraph{Discussion of these results in the context of using semiclassical analysis in the numerical analysis of the Helmholtz equation.}

In the last 10 years, there has been growing interest in using results about the $k$-explicit analysis of the Helmholtz equation 
from \emph{semiclassical analysis} (a branch of \emph{microlocal analysis}) to design and analyse numerical methods for the Helmholtz equation\footnote{A closely-related activity is the design and analysis of numerical methods for the Helmholtz equation based on proving \emph{new} results about the $k\tendi$ asymptotics of Helmholtz solutions for polygonal obstacles; see \cite{ChLa:07, HeLaMe:13, HeLaCh:14, ChHeLaTw:15, He:15}.
}. The activity has occurred in, broadly speaking, four different directions:

\ben
\item The use of the results of Melrose and Taylor \cite{MeTa:85} -- on the rigorous $k\tendi$ asymptotics of the solution of the Helmholtz equation in the exterior of a smooth convex obstacle with strictly positive curvature -- to design and analyse $k$-dependent approximation spaces for integral-equation formulations \cite{DoGrSm:07, GaHa:11, AsHu:14, EcHa:16, EcOz:17, Ec:18}, 
\item The use of the results of Melrose and Taylor \cite{MeTa:85}, along with the work of Ikawa \cite{Ik:88} on scattering from several convex obstacles, to analyse algorithms for multiple scattering problems \cite{EcRe:09, AnBoEcRe:10}.
\item The use of bounds on the Helmholtz solution operator (also known as \emph{resolvent estimates}) due to Vainberg \cite{Va:75} (using the propagation of singularities results of Melrose and Sj\"ostrand \cite{MeSj:82}) and Morawetz \cite{Mo:75} to prove bounds on both $\normAinv_{\LtGt}$ and the inf-sup constant of the domain-based variational formulation \cite{ChMo:08, Sp:14, BaSpWu:16, ChSpGiSm:17}, and also to analyse preconditioning strategies \cite{GaGrSp:15}.
\item The use of identities originally due to Morawetz \cite{Mo:75} to prove coercivity of $\opA$ \cite{SpKaSm:15} and to introduce new coercive formulations of Helmholtz problems \cite{SpChGrSm:11, MoSp:14, GaMo:17a, GaMo:17b, DiMoSp:18}.
\een
The results of the present paper arise from a fifth direction, namely
using estimates on the restriction of quasimodes 
of the Laplacian to hypersurfaces from 
\cite{Tat, BGT, T, HTacy, christianson2014exterior, T14} to prove sharp $k$-explicit bounds on $S_k, D_k$ and $D'_k$ as operators from $\LtG$ to $\HoG$.
We state these sharp $k$-explicit bounds in \S\ref{sec:L2H1} below, and they are proved in the companion paper \cite{GaSp:18}.
In the present paper, we use these new results, in conjunction with the results in Points 3 and 4 above, to obtain answers to Q1 and Q2.

\subsection{Formulation of the problem}\label{sec:form}

\subsubsection{Geometric definitions.}

Let $\Oi\subset \Rea^d,$ $d=2$ or $3,$ be a bounded Lipschitz open set, 
such that the open complement $\Oe:= \Rea^d \setminus \overline{\Oi}$ is connected.
Let $H^1_{\text{loc}}(\Omega_+)$ denote the set of functions $v$ such that $\chi v \in H^1(\Omega_+)$ for every $\chi \in C^\infty_{\rm comp}(\overline{\Omega_+}):=\{ \chi|_{\Omega_+}  : \,\chi \in C^\infty(\Rea^d) \text{ is compactly supported}\}$.
Let $\gamma^{+}$ denote the trace operator from $\Omega_{+}$ to $\bound$. 
Let $\bn$ be the outward-pointing unit normal vector to $\Oi$ (i.e.~$\bn$ points \emph{out} of $\Oi$ and \emph{in} to $\Oe$), and let $\partial_n^+$ 
denote the normal derivative trace operator from $\Oe$ to $\bound$ that satisfies $\dnpu = \bn \cdot \gamma^+(\gu)$ when $u\in H^2_{\text{loc}}(\Oe)$.
(We also call $\gamma^+ u$ the Dirichlet trace of $u$ and $\partial_n^+ u$ the Neumann trace.)

\begin{definition}[Star-shaped, and star-shaped with respect to a ball]

\

(i) $\Oi$ is \emph{star-shaped with respect to the point $\bx_0\in \Oi$} if, whenever $\bx \in \Oi$, the segment $[\bx_0,\bx]\subset \Oi$.

(ii) $\Oi$ is \emph{star-shaped with respect to the ball $B_{a}(\bx_0)$} if it is star-shaped with respect to every point in $B_{a}(\bx_0)$.

(iii) $\Oi$ is \emph{star-shaped with respect to a ball} if there exists $a>0$ and $\bx_0\in\Oi$ such that $\Oi$ is star-shaped with respect to the ball $B_{a}(\bx_0)$.
\end{definition}

\begin{definition}[Nontrapping]
We say that $\Oi\subset \Rea^d, \,d=2, 3$ is 
\emph{nontrapping} if $\bound$ is smooth ($C^\infty$) and,
given $R$ such that $\overline{\Oi}\subset B_R(\bze)$, there exists a $T(R)<\infty$ such that 
all the billiard trajectories (in the sense of Melrose--Sj{\"o}strand~\cite[Definition 7.20]{MeSj:82})
that start in $\Oe\cap B_R(\bze)$ at time zero leave $\Oe\cap B_R(\bze)$ by time $T(R)$.
\end{definition}

\begin{definition}[Smooth hypersurface]
We say that $\Gamma\subset \R^d$ is a {\em smooth hypersurface} if there exists $\Gammaext$ a compact embedded 
 smooth $d-1$ dimensional submanifold of $\R^d$, possibly with boundary,
 such that $\Gamma$ is an open subset of $\Gammaext$, with $\Gamma$ strictly away from $\partial \Gammaext$, and the boundary of $\Gamma$ can be written as a disjoint union
\beqs
\partial \Gamma=\left(\bigcup_{\ell=1}^n Y_\ell\right)\cup \Sigma,
\eeqs
where each $Y_\ell$ is an open, relatively compact, smooth embedded manifold of dimension $d-2$ in $\Gammaext$, $\Gamma$ lies locally on one side of $Y_\ell$, and  $\Sigma$ is closed set with $d-2$ measure $0$ and $\Sigma \subset \overline{\bigcup_{l=1}^nY_l}$. We then refer to the manifold $\Gammaext$ as an extension of $\Gamma$. 
\end{definition}
\noi For example, when $d=3$, the interior of a 2-d polygon is a smooth hypersurface, with $Y_i$ the edges and $\Sigma$ the set of corner points.
\begin{definition}[Curved]\label{def:curved}
We say a smooth hypersurface is \emph{curved} if there is a choice of normal so that the second fundamental form of the hypersurface is everywhere positive definite.
\end{definition}
\noi Recall that the principal curvatures are the eigenvalues of the matrix of the second fundamental form in an orthonormal basis of the tangent space, and thus ``curved" is equivalent to the principal curvatures being everywhere strictly positive (or everywhere strictly negative, depending on the choice of the normal).
\begin{definition}[Piecewise smooth]\label{def:psh}
We say that a hypersurface $\Gamma$ is \emph{piecewise smooth} if $\Gamma=\cup_{i=1}^N \overline{\Gamma}_i$ where $\Gamma_i$ are smooth hypersurfaces 
and $\Gamma_i\cap \Gamma_j=\emptyset.$
\end{definition}
\begin{definition}[Piecewise curved]\label{def:pc}
We say that a piecewise-smooth hypersurface $\Gamma$ is \emph{piecewise curved} if $\Gamma$ is as in Definition \ref{def:psh} and each $\Gamma_j$ is curved.
\end{definition}

\subsubsection{The boundary value problem and integral equation formulation}

\begin{definition}[Sound-soft scattering problem]\label{def:SSSP}
Given $k>0$ and an incident plane wave $u^I(\bx) = \exp(\ri k\bx\cdot\ba)$ for some $\ba\in \Rea^d$ with $|\ba|= 1$, find $u^S \in C^2(\Omega_+) \cap H^1_{\text{\emph{loc}}}(\Omega_+)$ such that the total field $u:= u^I + u^S$ satisfies the Helmholtz equation \eqref{eq:Helm} in $\Oe$, $\gpu = 0$ on $\bound$, and $u^S$ satisfies the Sommerfeld radiation condition 
\beqs
\pdiff{u^S}{r}(\bx) - \ri k \,u^S(\bx) = 
o\bigg(\frac{1}{r^{(d-1)/2}}\bigg)
\eeqs
as $r := |\bx| \tendi$, uniformly in 
$\bx/r$.
\end{definition}

\noi The incident field in the sound-soft scattering problem of Definition \ref{def:SSSP} is a plane wave, but this could be replaced by a point source or, more generally, a solution of the Helmholtz equation in a neighbourhood of $\Omega$; see \cite[Definition 2.11]{ChGrLaSp:12}.

\paragraph{Obtaining the direct integral equation \eqref{eq:direct}.}
If $u$ satisfies the sound-soft scattering problem of Definition \ref{def:SSSP} then Green's integral representation implies that
\beq\label{eq:IR}
u(\bx)= u^I(\bx) - \int_\bound \Phi_k(\bx,\by)\dnpu(\by)\,\rd s(\by), \quad \bx \in \Oe,
\eeq
(see, e.g., \cite[Theorem 2.43]{ChGrLaSp:12}), where $\Phi_k(\bx,\by)$ is the fundamental solution of the Helmholtz equation given by 
\beqs
\Phi_k(\bx,\by)=\displaystyle\frac{\ri}{4}H_0^{(1)}\big(k|\bx-\by|\big), \,\,d=2,\quad\quad \Phi_k(\bx,\by) = \frac{\re^{\ri k |\bx-\by|}}{4\pi |\bx-\by|}, \,\,d=3
\eeqs
(note that we have chosen the sign of $\Phi_k(\bx,\by)$ so that $-(\Delta +k^2)\Phi_k(\bx,\by)= \delta(\bx-\by)$).
Taking the exterior Dirichlet and Neumann traces of \eqref{eq:IR} on $\bound$ and using the jump relations for the single- and double-layer potentials (see, e.g., \cite[Equation 2.41]{ChGrLaSp:12}) we obtain the integral equations
\beq\label{eq:BIE1}
S_k \dnpu = \gamma^+ u^I \quad\tand \quad \left( \half I + D_k'\right) \dnpu = \partial_n^+ u^I,
\eeq
where $S_k$ and $D'_k$ are the single-  and adjoint-double-layer operators defined by 
\begin{align}
S_k \phi(\bx) := &\int_\bound \Phi_k(\bx,\by) \phi(\by)\,\rd s(\by), \quad
D'_k \phi(\bx) := \int_\bound \frac{\partial \Phi_k(\bx,\by)}{\partial n(\bx)}  \phi(\by)\,\rd s(\by),
 \label{eq:SD'}
\end{align}
 for $\phi\in\LtG$ and $x\in\bound$. Later we will also need the definition of the double-layer potential,
\beq \label{eq:D}
D_k \phi(\bx) := \int_\bound \frac{\partial \Phi_k(\bx,\by)}{\partial n(\by)}  \phi(\by)\,\rd s(\by) \quad\tfor \phi\in\LtG \,\tand \, \bx\in\bound.
\eeq

The first equation in \eqref{eq:BIE1} is not
uniquely solvable when $-k^2$ is a Dirichlet eigenvalue of the
Laplacian in $\Oi$, and the second equation in \eqref{eq:BIE1} is not uniquely solvable
when $-k^2$ is a Neumann eigenvalue of the Laplacian in $\Oi$ (see, e.g., \cite[Theorem 2.25]{ChGrLaSp:12}).
The standard way to resolve this difficulty is to take a linear combination of the two equations, which yields the  integral equation
\eqref{eq:direct}
where $\opA$ is defined by \eqref{eq:scpo},
\beq\label{eq:f}
f_{k,\eta}:= \dnpu^I- \ri \eta\, \gamma^+ u^I,
\eeq
and we use the notation that $v:=\dnpu$ (this makes denoting the Galerkin solution below easier, since we then have $v_h$ instead of $(\dnpu)_h$).

The space $\LtG$ is a natural space for the practical solution of second-kind integral equations since it is self-dual, and, for $\eta\in \Rea\setminus\{0\}$, $\opA$ is a bounded invertible operator from $\LtG$ to itself \cite[Theorem 2.27]{ChGrLaSp:12}.  Furthermore the right-hand side $f_{k,\eta}$ is in $\LtG$ (since $u^I\in C^\infty(\overline{\Oe})$) and thus we consider the equation \eqref{eq:direct} as an equation in $\LtG$.

\paragraph{The Galerkin method.}

Given a finite-dimensional approximation space 
$\cV_N \subset \LtG$, the Galerkin method for the integral equation \eqref{eq:direct} is
\beq\label{eq:Galerkin}
\text{find } v_N \in \cV_N \text{ such that }\, \big(\opA v_N, w_N\big)_{\LtG}=\big(f_{k,\eta},w_N\big)_{\LtG} \, \text{ for all } w_N\in\cV_N.
\eeq
Let $\cV_N= \sspan\{\phi_i: i=1,\ldots,N\}$, let $v_N \in \cV_N$ be equal to $\sum_{j=1}^N V_j \phi_j$, and define $\bv\in \Com^N$ by  $\bv := (V_j)_{j=1}^N$. 
Then, with $\bA_{ij}:= (\opA \phi_j,\phi_i)_{\LtG}$ and $\bff_i:= (f_{k,\eta},\phi_i)_{\LtG}$, the Galerkin method \eqref{eq:Galerkin} is equivalent to solving the linear system $\bA \bv = \bff$.

We consider the $h$--version of the Galerkin method, and we then denote $\cV_N$ and $v_N$ by $\cV_h$ and $v_h$ respectively.
The main results for Q1 and Q2 will be stated under the following assumption on $\cV_h$.

\begin{assumption}[Assumptions on $\cV_h$]\label{ass:Vh}
$\cV_h$ is a space of piecewise polynomials of degree $p$ for some fixed $p\geq 0$ on shape-regular meshes of diameter $h$, with $h$ decreasing to zero
 (see, e.g., \cite[Chapter 4]{SaSc:11} for specific realisations).
Furthermore

\noi (a) if $w\in H^1(\bound)$ then
\beq \label{eq:appth}
\min_{w_h \in \cV_h}\N{w - w_h}_{\LtG} \lesssim  h\N{w}_{\HoG}, 
\eeq
(b) 
\beq\label{eq:normequiv}
\N{w_h}^2_{\LtG} \sim h^{d-1} \N{\bw}_2^2,
\eeq
where $\|\cdot\|_2$ denotes the $l_2$ (i.e.~euclidean) vector norm.
\end{assumption}

\bre[For what situations is Assumption \ref{ass:Vh} proved?]
Part (a) is proved for subspaces consisting of piecewise-constant basis functions in \cite[Theorem 4.3.19]{SaSc:11} when $\bound$ is a polyhedron or curved (in the sense of Assumptions 4.3.17 and 4.3.18, respectively, in \cite{SaSc:11}) and in \cite[Theorem 10.4]{St:08} when $\Oi$ is a piecewise-smooth Lipschitz domain.
Part (a) is proved for subspaces consisting of continuous piecewise-polynomials of degree $p\geq 1$ (in the sense of \cite[Definition 4.1.36]{SaSc:11}) in \cite[Theorem 4.3.28]{SaSc:11}.

Part (b) is proved for subspaces consisting of piecewise-linear basis function in \cite[Lemma 10.5]{St:08} when $\bound$ is piecewise-smooth and Lipschitz, and for more general subspaces in \cite[Theorem 4.4.7]{SaSc:11}. 
\ere

\subsection{Statement of the main results and discussion}

The results concerning Q1 are stated in \S\ref{sec:r2}, and the results concerning Q2 are stated in \S\ref{sec:r3}.

\subsubsection{Results concerning Q1}\label{sec:r2}

\begin{theorem}[Sufficient conditions for the Galerkin method to be quasi-optimal]\label{thm:1}
Let $u$ be the solution of the sound-soft scattering problem of  Definition \ref{def:SSSP} and let $v:=\dnpu$.
Let $|\eta|\sim k$, and let $\cV_h$ satisfy Part (a) of Assumption \ref{ass:Vh}.

(a) 
If \emph{either} (i) $\Oi$ is nontrapping, \emph{or} (ii) $\Oi$ is star-shaped with respect to a ball
and $\bound$ is $C^{2,\alpha}$ and piecewise smooth, then given $k_0>0$, there exists a $C>0$ (independent of $k$ and $h$) such that if
\beq\label{eq:thres1}
h k^{3/2} \leq C, \quad d=2,\quad h k^{3/2}\log k \leq C, \quad d=3,
\eeq
then the Galerkin equations \eqref{eq:Galerkin}  have a unique solution which satisfies
\beq\label{eq:qothm}
\N{v-v_h}_{\LtG} \lesssim \min_{w_h\in \cV_h} \N{v-w_h}_{\LtG}
\eeq
for all $k\geq k_0$.

(b) In case (ii) above, if additionally $\bound$ is piecewise curved, then 
given $k_0>0$, there exists a $C>0$ (independent of $k$ and $h$) such that if
\beq\label{eq:thres2}
h k^{4/3}\log k \leq C, \quad d=2,3
\eeq
then \eqref{eq:qothm} holds.

(c) If $\Oi$ is convex and $\bound$ is $C^\infty$ and curved then 
given $k_0>0$, there exists a $C>0$ (independent of $k$ and $h$) such that if
\beq\label{eq:thres3}
h k^{4/3} \leq C, \quad d=2,3
\eeq
then \eqref{eq:qothm} holds.
\end{theorem}

Having established quasi-optimality, it is then natural to ask how the best approximation error $\min_{w_h\in \cV_h} \N{v-w_h}_{\LtG}$ depends on $k$, $h$, and $\|v\|_\LtG$.

\begin{theorem}[Bounds on the best approximation error]\label{thm:2}
Let $u$ be the solution of the sound-soft scattering problem of Definition \ref{def:SSSP} and let $v:=\dnpu$.
Let $\cV_h$ satisfy Assumption \ref{ass:Vh}. 

(a) If $\bound$ is $C^{2,\alpha}$ and piecewise smooth, then, given $k_0>0$, 
\beq\label{eq:bae2}
\min_{w_h \in \cV_h}\N{v-w_h}_{\LtG} \lesssim h A(k) \N{v}_{\LtG}
\eeq
with 
$A(k)= k^{5/4}\log k$,
for all $k\geq k_0$. 

(b) If $\bound$ is piecewise curved, then, 
given $k_0>0$, \eqref{eq:bae2} holds with
$A(k)= k^{7/6}\log k$, 
for all $k\geq k_0$. 

(c) If $\Oi$ is convex and $\bound$ is $C^\infty$ and curved, then,  given $k_0>0$,
\eqref{eq:bae2} holds with 
$A(k)= k$, 
for all $k\geq k_0$. 
\end{theorem}

Combining Theorems \ref{thm:1} and \ref{thm:2} we can obtain bounds on the relative error of the Galerkin method. For brevity, we only state the ones corresponding to cases (a) and (c) in Theorems \ref{thm:1} and \ref{thm:2}.

\begin{corollary}[Bound on the relative errors in the Galerkin method]\label{cor:1}
Let $u$ be the solution to the sound-soft scattering problem, let $|\eta|\sim k$, and let $\cV_h$ satisfy Part (a) of Assumption \ref{ass:Vh}.

(a) 
If \emph{either} (i) $\Oi$ is nontrapping, \emph{or} (ii) $\Oi$ is star-shaped with respect to a ball
and $\bound$ is $C^{2,\alpha}$ and piecewise smooth, then given $k_0>0$, there exists a $C>0$ (independent of $k$ and $h$) such that if $h$ and $k$ satisfy \eqref{eq:thres1}
then 
the Galerkin equations \eqref{eq:Galerkin} have a unique solution which satisfies
\beqs
\frac{\N{v-v_h}_{\LtG}}{\N{v}_{\LtG}} 
\lesssim 
\begin{cases}
k^{-1/4}\log k ,& d=2,\\
k^{-1/4}, & d=3,
\end{cases}
\eeqs
for all $k\geq k_0$.

(b) If $\Oi$ is convex and $\bound$ is $C^\infty$ and curved, then given $k_0>0$, there exists a $C>0$ (independent of $k$ and $h$) such that if
$h k^{4/3} \leq C$ 
the Galerkin equations \eqref{eq:Galerkin} have a unique solution which satisfies
\beqs
\frac{\N{v-v_h}_{\LtG}}{\N{v}_{\LtG}} 
\lesssim \frac{1}{k^{1/3}}
\eeqs
for all $k\geq k_0$.
\end{corollary}

\bre[The main ideas behind the proofs of Theorems \ref{thm:1} and \ref{thm:2}]\label{rem:ideaQ1}
The proof of Theorem \ref{thm:1} uses the classic projection-method analysis of second-kind integral equations (see, e.g., \cite[Chapter 3]{At:97}), with $\opA$ treated as a compact perturbation of $\half I$. 
In \cite{GrLoMeSp:15}, this argument was used to reduce the question of finding $k$-explicit bounds on the mesh threshold $h$ for $k$-independent quasi-optimality to hold to finding $k$-explicit bounds on 
\beqs
\|S_k\|_{\LtG\rightarrow\HoG}, \quad \|D'_k\|_{\LtG\rightarrow\HoG}, \quad\tand \quad\normAinv_{\LtGt}.
\eeqs
We use the new, sharp bounds on the first two of these norms from \cite{GaSp:18}, quoted here as Theorem \ref{thm:L2H1}, and the sharp bounds on the third of these norms from 
\cite[Theorem 1.13]{BaSpWu:16} (for nontrapping obstacles) and \cite[Theorem 4.3]{ChMo:08} (for obstacles that are star-shaped with respect to a ball).

The bounds of Theorem \ref{thm:2} are proved by showing that
\beq\label{eq:H1L2}
\N{v}_{\HoG}\lesssim A(k) \N{v}_{\LtG},
\eeq
and then using the approximation theory result \eqref{eq:appth}.
The bound \eqref{eq:H1L2} is obtained from the integral equation \eqref{eq:direct} using the second-kind-structure of the equation and 
the $\LtG\rightarrow \HoG$ bounds on $S_k$ and $D'_k$ from Theorem \ref{thm:L2H1}.
\ere

\bre[Comparison to previous results]
Theorems \ref{thm:1} and \ref{thm:2} and Corollary \ref{cor:1} sharpen previous results in \cite{GrLoMeSp:15}:~the mesh thresholds for quasi-optimality in Theorem \ref{thm:1} are sharper than the corresponding ones in \cite{GrLoMeSp:15}, and the results are valid for a wider class of obstacles.

This sharpening is due to the new, sharp bounds on  $\LtG\rightarrow\HoG$ norms of $S_k$, $D_k$, and $D'_k$ from \cite{GaSp:18},
and the widening of the class of obstacles is due to the bound on 
$\normAinv_{\LtGt}$ for nontrapping obstacles from \cite[Theorem 1.13]{BaSpWu:16}.
In more detail:
Theorem 1.4 of \cite{GrLoMeSp:15} is the analogue of our Theorem \ref{thm:1} except that the former is only valid when $\Oi$ is star-shaped with respect to a ball and $C^{2,\alpha}$ and the mesh threshold is $hk^{(d+1)/2}\leq C$. Comparing this result to Theorem \ref{thm:1} we see that we've sharpened the threshold in the $d=3$ case, expanded the class of obstacles to nontrapping ones, and added the additional results (b) and (c). 
Theorem \ref{thm:2} on the best approximation error is again proved using the $\LtG\rightarrow\HoG$-bounds from \cite{GaSp:18} and thus we see similar improvements over the corresponding theorem in \cite{GrLoMeSp:15} (\cite[Theorem 1.3]{GrLoMeSp:15}).

As discussed in Remark \ref{rem:ideaQ1}, both the present paper and \cite{GrLoMeSp:15} use the classic projection-method argument to obtain $k$-explicit results about quasi-optimality of the $h$-BEM. There are two other sets of results about quasi-optimality of the $h$-BEM in the literature:
\ben
\item[(a)] results that use coercivity \cite{DoGrSm:07, SpChGrSm:11, SpKaSm:15}, and
\item[(b)] results that give sufficient conditions for quasi-optimality to hold in terms of how well the spaces $\cV_h$ approximate the solution of certain adjoint problems \cite{BaSa:07, LoMe:11, Me:12}.
\een
These two sets of results are discussed in detail in \cite[pages 181--182]{GrLoMeSp:15} and \cite[\S4.2]{GrLoMeSp:15} respectively, and neither give results as strong as those in Theorem \ref{thm:1}.

Finally, in this paper we have only considered the $h$-BEM; a thorough $k$-explicit analysis of the $hp$-BEM for the exterior Dirichlet problem was conducted in \cite{LoMe:11} and \cite{Me:12}. In particular, this analysis, combined with the bound on $\normAinv_{\LtGt}$ for nontrapping obstacles from \cite[Theorem 1.13]{BaSpWu:16},
proves that $k$-independent quasi-optimality can be obtained for nontrapping obstacles through a choice of $h$ and $p$ that keeps the total number of degrees of freedom proportional to $k^{d-1}$ \cite[Corollaries 3.18 and 3.19]{LoMe:11}.
\ere

\bre[How sharp are the quasioptimality results?]
Numerical experiments in \cite[\S5]{GrLoMeSp:15} show that for a wide variety of obstacles (including certain mildly-trapping obstacles) the $h$-BEM is quasi-optimal with constant independent of $k$ (i.e.~\eqref{eq:qothm} holds), when $hk\sim 1$. The closest we can get to proving this is the result for strictly convex obstacles in Theorem \ref{thm:1} part (c), with the threshold being $hk^{4/3}\leq C$. 
The recent results of Baydoun and Marburg \cite{Ma:16, BaMa:17, Ma:17, BaMa:18}, however, give examples of cases where $hk\sim 1$ is not sufficient to keep the error bounded as $k\tendi$.
\ere

\subsubsection{Result concerning Q2}\label{sec:r3}

We now consider solving the linear system $\bA \bv = \bff$ with the generalised minimum residual method (GMRES) introduced by Saad and Schultz in \cite{SaSc:86}; for details of the implementation of this algorithm, see, e.g., \cite{Sa:03, Gr:97}.

\begin{theorem}[A bound on the number of GMRES iterations]\label{thm:GMRES1}
Let $\Oi$ be a 2- or 3-d convex obstacle whose boundary $\bound$ is analytic and curved. 
Let $\cV_h$ satisfy Part (b) of Assumption \ref{ass:Vh}, let the Galerkin matrix corresponding to \eqref{eq:Galerkin} be denoted by $\bA$, and consider GMRES applied to $\bA \bv= \bff$ 

There exist constants  $\eta_0>0$ and $k_0>0$ (with $\eta_0=1$ if $\Omega$ is a ball) such that if $k\geq k_0$ and $\eta_0 k \leq \eta\lesssim k$, then,
given $0<\eps <1$, there exists a $C$ (independent of $k$, $\eta$, and $\eps$) such that if
\beq\label{eq:boundonm}
m\geq C k^{1/3}\log\left( \frac{12}{\eps}\right),
\eeq
then the $m$th GMRES residual $\br_m:= \bA \bv_m -\bff$ satisfies
\beqs
\frac{\N{\br_m}_2}{\N{\br_0}_2}\leq \eps,
\eeqs
where $\|\cdot\|_2$ denotes the $l_2$ (i.e.~euclidean) vector norm.
\end{theorem}

\noi In other words, Theorem \ref{thm:GMRES1} states that, for convex, analytic, curved $\Omega$, the number of iterations growing like $k^{1/3}$ is a sufficient condition for GMRES to maintain accuracy as $k\tendi$.

\bre[How sharp is the result of Theorem \ref{thm:GMRES1}?]
Numerical experiments in \S\ref{sec:num} show that for the sphere 
the number of GMRES iterations grows like $k^{0.29}$, and for an ellipsoid they grow like $k^{0.28}$.
The bound in Theorem \ref{thm:GMRES1} is therefore effectively sharp (at least for the range of $k$ considered in the experiments).
\ere

\bre[The main ideas behind the proof of Theorem \ref{thm:GMRES1}]\label{rem:ideaQ2}
The two ideas behind Theorem \ref{thm:GMRES1} are that:

(a) A sufficient (but not necessary) condition for iterative methods to be well behaved is that the \emph{numerical range} (also known as the \emph{field of values}) of the matrix is bounded away from zero, and in this case the Elman estimate \cite{El:82, EiElSc:83} and its refinement due to 
Beckermann, Goreinov, and Tyrtyshnikov \cite{BeGoTy:06} can be used to bound the number of GMRES iterations in terms of (i) the distance of the numerical range to the origin, and (ii) the norm of the matrix.

(b) When $\Oi$ is convex, $C^3$, piecewise analytic, and $\bound$ is curved, \cite{SpKaSm:15} proved that $\opA$ is coercive for sufficiently large $k$ (with $\eta\sim k$). The $k$-dependence of the coercivity constant, along with the $k$-dependence of $\|\opA\|_\LtGt$ then give the information needed about the numerical range of the Galerkin matrix $\bA$ required in (a).
\ere

\bre[Comparison to previous results]
The bound $m\gtrsim k^{2/3}$ when $\bound$ is a sphere was stated in \cite[\S1.3]{SpKaSm:15}; this bound was obtained using the original Elman estimate (see Remark \ref{rem:Elman} below), and the fact that the sharp bound $\|\opA\|_{\LtGt}\lesssim k^{1/3}$ was known for the circle and sphere; see \cite[\S5.4]{ChGrLaSp:12}. To our knowledge, there are no other $k$-explicit bounds in the literature on the number of GMRES iterations required to achieve a prescribed accuracy for a Helmholtz BIE. 
The closest related work is \cite{ChRaRo:02}, which uses a second-kind integral equation to solve the Helmholtz equation in a half-plane with an impedance boundary condition. The special structure of this integral equation allows a two-grid iterative method to be used, and \cite{ChRaRo:02} proves that 
there exists $C>0$ such that if $k h\leq C$, then, after seven iterations, the difference between the solution and the Galerkin solution computed via the iterative method is bounded independently of $k$ and $h$.
\ere

\bre[Translating the results to the indirect equation \eqref{eq:indirect}]\label{rem:indirect}
Instead of using Green's integral representation \eqref{eq:IR} to formulate the sound-soft scattering problem as the  integral equation \eqref{eq:direct}, one can 
pose the ansatz that the scattered field satisfies
\beqs
u^S(\bx)=\int_\bound \frac{\partial \Phi_k(\bx,\by)}{\partial n(\by)}  \phi(\by)\,\rd s(\by) - \ri \eta \int_\bound \Phi_k(\bx,\by)\phi(\by)\,\rd s(\by)
\eeqs
for $\bx\in \Oe$, $\phi\in \LtG$, and $\eta\in \Rea\setminus \{0\}$. Imposing the boundary condition $\gamma^+u^S=-\gamma^+ u^I$ on $\bound$ and using the jump relations for the single- and double-layer potentials leads to the integral equation \eqref{eq:indirect} where $\opABW$ is defined by \eqref{eq:scpo} and $g= - \gamma^+ u^I$. 
One can show that $\opABW$ and $\opA$ are adjoint 
with respect to the real-valued $\LtG$ inner product (see, e.g., \cite[Equation 2.37, Remark 2.24, \S2.6]{ChGrLaSp:12}), and so their norms are equal, the norms of their inverses are equal, and if one is coercive then so is the other (with the same coercivity constant).
These facts imply that the results of Theorems \ref{thm:1} and Theorem \ref{thm:GMRES1} hold for the indirect equation \eqref{eq:indirect}.

The bounds on the best approximation error in Theorem \ref{thm:2} hold for the indirect equation \eqref{eq:indirect} with (a) $A(k) = k^{3/2}$ for $d=2$, $A(k)=k^{3/2}\log k$ for $d=3$, (b) $A(k)=k^{4/3}\log k$, and (c) $A(k)=k^{4/3}$. These powers of $k$ are all slightly higher than those for the direct equation; the reason for this is that we have more information about the unknown in the direct equation (since it is $\dnpu$) 
than about the unknown $\phi$ in the indirect equation. Indeed, one can express $\phi$ in terms of the difference of solutions to interior and exterior boundary value problems -- see \cite[Page 132]{ChGrLaSp:12}  -- but it is harder to make use of this fact than for the direct equation.
\ere

\bre[Translating the results to the general exterior Dirichlet problem]
The results of Theorems \ref{thm:1} and \ref{thm:GMRES1} are independent of the right-hand side of the integral equation \eqref{eq:direct}, and therefore hold for the general Dirichlet problem with Dirichlet data in $H^1(\bound)$ (this assumption is needed so that $\opA$ can still be considered as an operator on $\LtG$; see, e.g., \cite[\S2.6]{ChGrLaSp:12}).
The results of Theorem \ref{thm:2} and Corollary \ref{cor:1}, however, do not immediately hold for the general Dirichlet problem, since they use the particular form of the right-hand side in \eqref{eq:f}.
\ere

\paragraph{Outline of the paper}
In \S\ref{sec:L2H1} we recap the sharp $\LtG\rightarrow \HoG$ bounds from the companion paper \cite{GaSp:18}.
In \S\ref{sec:Q1} we prove Theorems \ref{thm:1} and \ref{thm:2} (the results concerning Q1).
In \S\ref{sec:Q2} we prove Theorem \ref{thm:GMRES1} (the result concerning Q2), and then in \S\ref{sec:num} we give numerical experiments showing that Theorem \ref{thm:GMRES1} is effectively sharp in its $k$-dependence.

\section{Recap of the $\LtG\rightarrow \HoG$ bounds from \cite{GaSp:18}}\label{sec:L2H1}

The following result was proved in \cite[Theorem 2.10]{GaSp:18}. 
In stating this result, we use the weighted $\HoG$ norm
\beq\label{eq:HohG}
\N{w}_{\HohG}^2:= k^{-2}\N{\nT w}_{\LtG}^2+ \N{w}_{\LtG}^2,
\eeq
in contrast to the usual $\HoG$ norm 
\beqs
\N{w}_{\HoG}^2:=\N{\nT w}_{\LtG}^2+ \N{w}_{\LtG}^2,
\eeqs
where $\nT$ is the surface gradient operator on $\bound$ (see, e.g., \cite[Page 276]{ChGrLaSp:12}); we note that the use of such weighted norms is standard in both the semiclassical and numerical analysis of the Helmholtz equation, and reflects the fact that we expect to incur a power of $k$ every time we take a derivative of a Helmholtz solution; see, e.g., \cite[Remark 3.8]{MoSp:14}.

\begin{theorem}[Bounds on $\|S_k\|_{\LtG\rightarrow\HohG}$, $\|D_k\|_{\LtG\rightarrow\HohG}$, $\|D'_k\|_{\LtG\rightarrow\HohG}$]
\label{thm:L2H1}

\

\noi Let $\Omega$ be a bounded Lipschitz open set such that the open complement $\Oe:= \Rea^d\setminus \overline{\Omega}$ is connected.

\noi (a) 
If $\bound$ is a piecewise-smooth hypersurface (in the sense of Definition \ref{def:psh}), then, 
given $k_0>1$, there exists $C>0$ (independent of $k$) such that
\begin{equation}
\label{eqn:optimalFlatSl_new}
\|\Sk\|_{L^2(\bound)\to  H_{k}^1(\bound)}  \leq 
C\,k^{-1/2}\,\log k.
\end{equation}
for all $k\geq k_0$.
Moreover, if $\bound$ is piecewise curved (in the sense of Definition \ref{def:pc}), 
 then, given $k_0>1$, the following stronger estimate holds for all $k\geq k_0$
\begin{equation}\label{eqn:optimalConvexSl_new}
\N{S_k}_{\LtG\rightarrow{H}_k^1(\bound)} \leq C k^{-2/3}\log k.
\end{equation} 

\noi (b) If $\bound$ is a piecewise smooth, $\Dkreg$ hypersurface, for some $\alpha>0$, then, given $k_0>1$, there exists $C>0$ (independent of $k$) such that
\begin{equation*}
\N{D_k}_{\LtG\rightarrow{H}_k^1(\bound)} + \N{D'_k}_{\LtG\rightarrow{H}_k^1(\bound)} \leq C k^{1/4}\log k.
\end{equation*}
Moreover, if $\bound$ is piecewise curved, 
then, given $k_0>1$, there exists $C>0$ (independent of $k$) such that the following stronger estimates hold for all $k\geq k_0$
\begin{equation*}
\N{D_k}_{\LtG\rightarrow H^1_k(\bound)} + \N{D'_k}_{\LtG\rightarrow H^1_k(\bound)} \lesssim k^{1/6}\log k.
\end{equation*}

\noi (c) 
If $\Oi$ is convex and $\bound$ is $C^\infty$ and curved (in the sense of Definition \ref{def:curved}) then, given $k_0>1$, there exists $C$ such that, for $k\geq k_0$,
\begin{align*}
&\qquad\qquad\N{S_k}_{\LtG\rightarrow H_k^1(\bound)} \leq C k^{-2/3}, \\
&\N{D_k}_{\LtG\rightarrow H_k^1(\bound)}+\N{D'_k}_{\LtG\rightarrow H_k^1(\bound)} \leq C.
\end{align*}
\end{theorem}

The requirement in Part (b) of Theorem \ref{thm:L2H1} that $\bound$ is $C^{2,\alpha}$ arises since this is the regularity required of $\bound$ for $D_k$ and $D'_k$ to map $\LtG$ to $\HoG$; see \cite[Theorem 4.2]{Ki:89},  \cite[Theorem 3.6]{CoKr:98}.

The bounds in Theorem \ref{thm:L2H1} contain $k$-explicit $\LtG\rightarrow\LtG$ bounds on $S_k, D_k$ and $D'_k$. These 
$\LtG\rightarrow\LtG$ bounds were originally proved in 
\cite[Appendix A]{HaTa:15} and \cite{Ga:15} (and the realisation that these $\LtG\rightarrow\LtG$ bounds could be extended to  $\LtG\rightarrow\HoG$ bounds was the motivation for \cite{GaSp:18}).

\bre[Sharpness of the bounds in Theorem \ref{thm:L2H1}]
In  \cite[\S3]{GaSp:18} it is shown that, modulo the factor $\log k$, all of the bounds in Theorem \ref{thm:L2H1} are sharp (i.e.~the powers of $k$ in the bounds are optimal). 
The sharpness (modulo the factor $\log k$) of the $\LtG\rightarrow\LtG$ bounds contained in Theorem \ref{thm:L2H1} was proved in 
\cite[\S A.2-A.3]{HaTa:15}. Earlier work in \cite[\S4]{ChGrLaLi:09} proved the sharpness of some of the $\LtG\rightarrow\LtG$ bounds in 2-d; we highlight that \cite[\S3]{GaSp:18} and \cite[\S A.2-A.3]{HaTa:15} contain the appropriate generalisations to multidimensions of some of the arguments of \cite[\S4]{ChGrLaLi:09} (in particular \cite[Theorems 4.2 and 4.4]{ChGrLaLi:09}).
\ere

\bre[Sharp bounds on $S_k$ when $d=2$]
When $d=2$ and $\bound$ is Lipschitz, the sharp bound
\beq\label{eq:SbL2}
\|S_k\|_{\LtGt}\lesssim k^{-1/2}
\eeq
was proved using the Riesz--Thorin interpolation theorem in \cite[Theorem 3.3]{ChGrLaLi:09} 
and by the Schur test in \cite[Theorem 6]{GaSm:15}. Similarly,  the sharp bound 
\beq\label{eq:Sb}
\N{S_k}_{\LtGtH} \lesssim k^{1/2}
\eeq
was proved using the Riesz--Thorin interpolation theorem in \cite[Theorem 1.6]{GrLoMeSp:15}.
\ere

\section{Proofs of Theorems \ref{thm:1}, \ref{thm:2} (the results concerning Q1)}\label{sec:Q1}

\subsection{Proof of Theorem \ref{thm:1}}

The heart of the proof of Theorem \ref{thm:1} is the following lemma.

\begin{lemma}\label{lem:qo1}
There exists a $\tilde{C}>0$ such that under the condition 
\beq\label{eq:thm1stara}
h \N{D_k'-\ri \eta S_k}_{\LtG\rightarrow \HoG} \normAinv_{\LtGt} \leq \widetilde{C}
\eeq
the Galerkin equations \eqref{eq:Galerkin} have a unique solution satisfying \eqref{eq:qothm}.
\end{lemma}

The presence of $\normAinv_{\LtGt}$ in \eqref{eq:thm1stara} means that before proving Theorem \ref{thm:1} using Lemma \ref{lem:qo1} we need to recall the following bounds on $\normAinv_{\LtGt}$.

\begin{theorem}\textbf{\emph{(\cite[Theorem 4.3]{ChMo:08}, \cite[Theorem 1.13]{BaSpWu:16})}}\label{thm:Ainv}
If $|\eta|\sim k$ and \emph{either} $\Oi$ is star-shaped with respect to a ball and $C^2$ in a neighbourhood of almost every point on $\Gamma$  \emph{or} $\Oi$ is nontrapping, then, given $k_0>0$, $\normAinv_\LtGt\lesssim 1$ for all $k\geq k_0$.
\end{theorem}

\bpf[Proof of Theorem \ref{thm:1} using Lemma \ref{lem:qo1}]
Using the triangle inequality, a sufficient condition for \eqref{eq:thm1stara} to hold is 
\beq\label{eq:thm1star}
h\left( \N{D_k'}_{\LtG\rightarrow \HoG} + |\eta|\N{S_k}_{\LtG\rightarrow \HoG}\right) \normAinv_{\LtGt} \leq \widetilde{C}.
\eeq
In \cite[Remark 2.22]{GaSp:18} it is shown that the $\LtG\rightarrow \HoG$ norms of $D_k'$ and $S_k$ are maximised in different regions of phase space, and thus we do not lose anything by using the triangle inequality, i.e.,~\eqref{eq:thm1star} is no less sharp than \eqref{eq:thm1stara} in terms of $k$-dependence.

The mesh thresholds \eqref{eq:thres1}, \eqref{eq:thres2}, \eqref{eq:thres3} then follow from using the bound $\normAinv_\LtGt\lesssim 1$ from Theorem \ref{thm:Ainv} and 
the different bounds on  $\|D_k'\|_{\LtG\rightarrow\HoG}$ and $\|S_k\|_{\LtG\rightarrow\HoG}$ in Theorem \ref{thm:L2H1} (recalling the definition of the $H^1_k(\bound)$ norm in \eqref{eq:HohG}),
apart from when $d=2$ when we use the bound on $S_k$ \eqref{eq:Sb} instead of \eqref{eqn:optimalFlatSl_new}.
\epf

\

To prove Theorem \ref{thm:1} we therefore only need to prove Lemma \ref{lem:qo1}. This was proved in \cite[Corollary 4.1]{GrLoMeSp:15}, but since the proof is short we repeat it here for completeness.

We first introduce some notation: let $P_h$ denote the orthogonal projection from $\LtG$ onto $\cV_h$ (see, e.g, \cite[\S3.1.2]{At:97}); then the Galerkin equations \eqref{eq:Galerkin} are equivalent to the operator equation
\beq\label{eq:Galerkin_Ph}
P_h \opA v_h = P_h f_{k,\eta}.
\eeq
The proof requires us to treat $\opA$ as a (compact) perturbation of the identity, and thus we let $L_{k,\eta}:= D_k' - \ri \eta S_k$. Furthermore, to make the notation more concise, we let $\lambda=1/2$. Therefore, the left-hand side of \eqref{eq:Galerkin_Ph} becomes $(\lambda I + P_h L_{k,\eta})v_h$, and the question of existence of a solution to \eqref{eq:Galerkin_Ph} boils down to the invertibility of $(\lambda I + P_h L_{k,\eta})$. Note also that, using the $P_h$ notation, the best approximation error on the right-hand side of \eqref{eq:qothm} is $\|(I-P_h)v\|_{\LtG}$.

The heart of the proof of Lemma \ref{lem:qo1} is the following lemma.
\begin{lemma}\label{lem:qo2}
If
\begin{equation}
\label{eq:6abscond} 
\N{(I - P_h) L_{k,\eta}}_{\LtGt} \normAinv_\LtGt
\leq   \frac{\delta }{1+ \delta}
\end{equation}
for some $\delta>0$, then 
the Galerkin equations have a unique solution, $v_h$, which satisfies the quasi-optimal error estimate 
\begin{equation}
\N{v - v_h}_{\LtG} \leq  \lambda  (1+ \delta)   
\normAinv_\LtGt
\|(I-P_h)v\|_{\LtG}.
\label{eqn:6lquasiopt} 
 \end{equation}
\end{lemma}

\begin{proof}[Proof of Lemma \ref{lem:qo1} using Lemma \ref{lem:qo2}]
By the polynomial-approximation result \eqref{eq:appth},
\beqs
\N{(I-P_h)L_{k,\eta}}_{\LtGt} \lesssim h \N{L_{k,\eta}}_{\LtG\rightarrow \HoG}.
\eeqs
Therefore, choosing, say, $\delta=1$, we find that there exists a $\tilde{C}>0$ such that \eqref{eq:thm1stara} implies that \eqref{eq:6abscond} holds.
\end{proof}

\

\noi Thus, to prove Theorem \ref{thm:1}, we only need to prove Lemma \ref{lem:qo2}.

\

\begin{proof}[Proof of Lemma \ref{lem:qo2}]
Since
\begin{align*}
\lambda I + P_h L_{k,\eta} &= \lambda I + L_{k,\eta} - (I-P_h)L_{k,\eta}
= \left( \lambda I + L_{k,\eta}\right) \left( I - \left(\lambda I+ L_{k,\eta}\right)^{-1} (I-P_h)L_{k,\eta}\right),
\end{align*}
if
\beqs
\big\|\left(\lambda I+L_{k,\eta}\right)^{-1} (I-P_h)L_{k,\eta}\big\|_{\LtGt}<1,
\eeqs
then $(\lambda I + P_h L_{k,\eta})$ is invertible using the classical result that $I-A$ is invertible if $\N{A}<1$. In this abstract setting $\|(I-A)^{-1}\|\leq (1-\N{A})^{-1}$, and thus if 
\eqref{eq:6abscond} holds we have
\begin{align}\nonumber
\big\|(\lambda I + P_h L_{k,\eta})^{-1}\big\|_{\LtGt} &\leq \big\|\left(\lambda I +L_{k,\eta}\right)^{-1}\big\|_{\LtGt} \frac{1}{1- \delta/(1+\delta) },\\
&= (1 + \delta) \, \big\|(\lambda I + L_{k,\eta})^{-1}\big\|_{\LtGt}.\label{eq:ivan1}
\end{align}
Writing the direct equation as $(\lambda I + L_{k,\eta})v=f$ and the Galerkin equation as $(\lambda I + P_h L_{k,\eta})v_h=P_h f$, we have
\begin{align}\nonumber 
v - v_h = v - (\lambda I + P_h L_{k,\eta})^{-1} P_h f &=  (\lambda I + P_h L_{k,\eta})^{-1}(\lambda v - P_h(f - L_kv)) \\&= \lambda \left(\lambda I+ P_h L_{k,\eta}\right)^{-1} (I- P_h )   v, \label{eq:igg1}.
  \end{align} 
and the result \eqref{eqn:6lquasiopt}  follows from using the bound \eqref{eq:ivan1} in \eqref{eq:igg1}.
\end{proof} 

\bre[Is there a better choice of $\eta$ than $|\eta|\sim k$?]
Theorem \ref{thm:1} is proved under the assumption that $|\eta|\sim k$. 
This choice of $\eta$ is widely recommended from studies of the condition number of $\opA$; see \cite[Chapter 5]{ChGrLaSp:12} for an overview of these.
From \eqref{eq:thm1star} we see that the best choice of $\eta$, from the point of view of obtaining the least-restrictive threshold for $k$-independent quasi-optimality, will minimise the $k$-dependence of
\beqs
\left( \N{D_k'}_{\LtG\rightarrow \HoG} + |\eta|\N{S_k}_{\LtG\rightarrow \HoG}\right) \normAinv_{\LtGt}.
\eeqs
There does not yet exist a rigorous proof that $|\eta|\sim k$ minimises this quantity, but 
\cite[\S7.1]{BaSpWu:16} outlines exactly the necessary results still to prove. 
\ere

\subsection{Proof of Theorem \ref{thm:2}}

\bpf[Proof of Theorem \ref{thm:2}]
By the polynomial-approximation result \eqref{eq:appth}, we only need to prove that the bound \eqref{eq:H1L2} hold with the different functions $A(k)$.
The idea is to take the $H^1$ norm of the integral equation \eqref{eq:direct} and then use the $\LtGt$ and $\LtG\rightarrow \HoG$ bounds contained in Theorems \ref{thm:L2H1}.

Taking the $H^1$ norm of \eqref{eq:direct} and using the notation that $\opA= \half I + L_{k,\eta}$ and $v:=\dnpu$, as in the proof of Theorem \ref{thm:1} above, we have that
\beqs
\N{v}_{\HoG} \lesssim \N{L_{k,\eta}}_{\LtG\rightarrow \HoG}\N{v}_{\LtG} + \N{f_{k,\eta}}_{\HoG}.
\eeqs
In this inequality, $\eta$ is just a parameter that appears in $L_{k,\eta}$ and $f_{k,\eta}$, with the equation holding for all values of $\eta$; in other words, the unknown $v(=\dnpu)$ does not depend on the value of $\eta$. We now seek to minimise the $k$-dependence of 
$\|L_{k,\eta}\|_{\LtG\rightarrow \HoG}$. Looking at the $k$-dependence of the $\LtG\rightarrow\HoG$-bounds on $S_k$ and $D'_k$ in Theorem \ref{thm:L2H1}, we see that, under each of the different geometric set-ups, the best choice is $\eta=0$, and thus 
\beq\label{eq:thm21}
\N{v}_{\HoG} \lesssim \N{D'_k}_{\LtG\rightarrow \HoG}\N{v}_{\LtG} + k^2
\eeq
where we have explicitly worked out the $k$-dependence of $\|f_{k,\eta}\|_{\HoG}$ using the definition \eqref{eq:f}.

Taking the $L^2$ norm of \eqref{eq:direct} (with $\eta=0$), 
and noting that $\|f_{k,\eta}\|_{\LtG}\sim k$, we have that
\beq\label{eq:thm22}
\big(1 + \N{D'_k}_{\LtGt}\big) \N{v}_{\LtG} \gtrsim k.
\eeq
Using \eqref{eq:thm22} in \eqref{eq:thm21}, we have 
\beq\label{eq:thm22a}
\N{v}_{\HoG} \lesssim \left(\N{D'_k}_{\LtG\rightarrow \HoG}+ k\big(1 + \N{D'_k}_{\LtGt}\big)\right) \N{v}_{\LtG}.
\eeq
Since the bounds on the $\LtG\rightarrow\HoG$-norm of $D'_k$ in Theorem \ref{thm:L2H1} are one power of $k$ higher that the $\LtG\rightarrow\LtG$-bounds in Theorem \ref{thm:L2H1}, using these norm bounds in \eqref{eq:thm22a} results in the bound $\|v\|_{\HoG}\lesssim A(k)\|v\|_{\LtG}$ with the functions of $A(k)$ as in the statement of theorem (and equal to
 the right-hand sides of the bounds on $\|D_k'\|_{\LtG\rightarrow \HoG}$ in Theorem \ref{thm:L2H1}).
\epf

\section{Proofs of Theorem \ref{thm:GMRES1} (the result concerning Q2)}\label{sec:Q2}

To prove Theorem \ref{thm:GMRES1} we need to recall (i) the result about coercivity of $\opA$ when $\Oi$ is convex, $C^3$, piecewise analytic, and curved from \cite{SpKaSm:15}, and (ii) the refinement of the Elman estimate in \cite{BeGoTy:06}.

\begin{theorem}[Coercivity of $\opA$ for $\Oi$ convex, $C^3$, piecewise analytic, and curved \cite{SpKaSm:15}]
\label{thm:coer} 
Let $\Oi$ be a convex domain in either 2- or 3-d whose boundary, $\bound$, is curved and is both $C^3$ and piecewise analytic.
Then there exist  constants $\eta_0>0$, $k_0>0$ (with $\eta_0=1$ when $\Omega$ is a ball) and a function of $k$, $\alpha_{k}>0$, such that for $k\geq k_0$ and $\eta\geq\eta_0 k$,
\beq\label{eq:coer}
\big|\big(\opA \phi, \phi\big)_{\LtG} \big|\geq 
\alpha_{k}\Vert \phi\Vert^2_{\LtG}\quad \tfa \phi\in\LtG,
\eeq
where
\beq\label{eq:coer_asym}
\alpha_{k} = \half - \cO\big( k^{-2/3}\log k\big) \quad\tas k \tendi.
\eeq
\end{theorem}

In stating this result we have used the bound on $S_k$
\eqref{eqn:optimalConvexSl_new}
 and \cite[Remark 3.3]{SpKaSm:15} to get the asymptotics \eqref{eq:coer_asym}. The fact that $\eta_0=1$ when $\Oi$ is a ball follows from \cite[Corollary 4.8]{SpChGrSm:11}.

\begin{theorem}[Refinement of the Elman estimate \cite{BeGoTy:06}]\label{thm:Elman}
Let $\bA$ be a matrix with $0\notin W(\bA)$, where 
$W(\bA):= \big\{ (\bA \bv, \bv) : \bv \in \Com^N, \|\bv\|_2=1\big\}$ is the \emph{numerical range} of $\bA$. 
Let $\beta\in [0,\pi/2)$ be defined such that
\beqs
\cos \beta = \frac{\mathrm{dist}\big(0, W(\bA)\big)}{\| \bA\|_{2}},
\eeqs
and let $\gamma_\beta$ be defined by 
\beq\label{eq:gamma_beta}
\gamma_\beta:= 2 \sin \left( \frac{\beta}{4-2\beta/\pi}\right).
\eeq
Suppose the matrix equation $\bA \bv = \bff$ is solved using GMRES, and let $\br_m:= \bA \bv_m - \bff$ be the $m$-th GMRES residual.
Then
\beq\label{eq:Elman2}
\frac{\|\br_m\|_{2}}{\|\br_0\|_{2}} \leq \left(2 + \frac{2}{\sqrt{3}}\right)\big(2+ \gamma_\beta\big) \,\gamma_\beta^m.
\eeq
\end{theorem}

When we apply the estimate \eqref{eq:Elman2} to $\bA$, we find that $\beta= \pi/2-\delta$, where $\delta= \delta(k)$ is such that $\delta\tendo$ as $k\tendi$.
We therefore specialise the result \eqref{eq:Elman2} to this particular situation in the following corollary.

\begin{corollary}\label{cor:Elman}
If $\beta= \pi/2-\delta$ with 
$0<\delta<\delta_0$, then there exists $C_1>0$ and $\delta_1>0$ (both independent of $\delta$) such that, for $0<\eps<1$, 
\beq\label{eq:corE1}
\text{if} \quad m\geq \frac{C_1}{\delta}\log\left(\frac{12}{\eps}\right)\quad\text{ then }  \quad\frac{\|\br_m\|_{D}}{\|\br_0\|_{D}}\leq \eps
\eeq
for all $0<\delta<\delta_1$.
\end{corollary}

\noi That is, choosing $m\gtrsim \delta^{-1}$ is sufficient for GMRES to converge in an $\delta$-independent way as $\delta\tendo$.

\

\bpf[Proof of Corollary \ref{cor:Elman}]
If $\beta= \pi/2-\delta$, with $\delta\rightarrow 0$, then $\cos\beta=\sin\delta = \delta+ \cO(\delta^3)$ as $\delta \tendo$. From the definition of the convergence factor  $\gamma_\beta$, \eqref{eq:gamma_beta}, we have
\beq\label{eq:factor}
\gamma_\beta:=2 \sin \left( \frac{\beta}{4-2\beta/\pi}\right)= 2 \sin \left( \frac{\pi}{6} -\frac{4\delta}{9}+\cO(\delta^2)\right)= 1 -\frac{4\delta}{3\sqrt{3}} + \cO(\delta^2)\quad\tas\,\, \delta\tendo,
\eeq
and then
\beqs
\log \gamma_\beta = -\frac{4\delta}{3\sqrt{3}} + \cO(\delta^2) \quad\tas\,\, \delta\tendo,
\eeqs
and so there exist  $C_2>0$ and $\delta_1>0$ such that
\beqs
\gamma_\beta^m = \re^{m\log\gamma_\beta}\leq \re^{-m\delta/C_2} \quad\tfa0< \delta\leq \delta_1.
\eeqs
The bound \eqref{eq:corE1} then follows from \eqref{eq:Elman2} since $(2 + 2/\sqrt{3})(2+ \gamma_\beta)<3(2 + 2/\sqrt{3})<12$.
\epf

\bre[Comparison of \eqref{eq:Elman2} with the original Elman estimate]\label{rem:Elman}
The estimate
\beq\label{eq:Elman}
\frac{\|\br_m\|_{2}}{\|\br_0\|_{2}} \leq \sin^m \beta
\eeq
was essentially proved in \cite{El:82, EiElSc:83} (see also
the review \cite[\S6]{SiSz:07} and the references therein).
When $\beta= \pi/2-\delta$, the convergence factor in \eqref{eq:Elman} is
\beqs
\sin\beta = \cos\delta 
= 1- \frac{\delta^2}{2}+\cO(\delta^4);
\eeqs
by comparing this to \eqref{eq:factor} we can see that \eqref{eq:Elman} is indeed a weaker bound.
\ere

\bpf[Proof of Theorem \ref{thm:GMRES1}]
The set up of the Galerkin method in \S\ref{sec:form} implies that, for any $v_N, w_N\in \cV_N$,
$( \opA v_N, w_N)_{\LtG} = (\bA \bv,\bw)_2$, where $(\cdot,\cdot)_2$ denotes the euclidean inner product on $l^2$. Therefore, the continuity of $\opA$ and the norm equivalence \eqref{eq:normequiv} implies that
\beq\label{eq:GMRES1}
|(\bA \bv, \bw)_2| \lesssim \N{\opA}_{\LtGt} h^{d-1} \N{\bv}_2 \N{\bw}_2 \quad\tfa \bv, \bw \in \Com^N.
\eeq
Furthermore, if $\opA$ is coercive with coercivity constant $\alpha_{k,\eta}$, i.e.,~\eqref{eq:coer} holds, 
then
\beq\label{eq:GMRES2}
|(\bA \bv, \bv)_2| \gtrsim \alpha_{k,\eta}h^{d-1} \N{\bv}^2_2  \quad\tfa \bv \in \Com^N.
\eeq
The bounds \eqref{eq:GMRES1} and \eqref{eq:GMRES2} together imply that the ratio $\cos \beta$ in \eqref{eq:Elman} satisfies
\beqs
\cos \beta \gtrsim \frac{\alpha_{k,\eta}}{\|\opA\|_\LtGt}.
\eeqs

Since $\Oi$ is $C^\infty$ and curved, the bound 
$\|\opA\|_{\LtGt}\lesssim k^{1/3}$ follows from the bounds in Theorem \ref{thm:L2H1} (recalling that $\eta_0 k\leq \eta\lesssim k$). Since $\bound$ is piecewise analytic, $C^3$, and curved, from Theorem \ref{thm:coer} there exists a $k_0>0$ such that $\alpha_{k,\eta}\sim 1$ for all $k\geq k_0$. Combining these two bounds we have $\cos\beta \gtrsim k^{-1/3}$ for all $k\geq k_0$
 and thus Corollary \ref{cor:Elman} holds with $\delta\sim  k^{-1/3}$ for all $k\geq k_0$; the result \eqref{eq:boundonm} then follows from \eqref{eq:corE1}.

Note that the assumption in the theorem that $\bound$ is analytic comes from the fact that if $\bound$ is both piecewise analytic and $C^\infty$, then $\bound$ must be analytic, where the notion of piecewise analyticity in Theorem \ref{thm:coer} is inherited from 
\cite[Definition 4.1]{ChSo:04}.
\epf 

\bre[The star-combined operator]\label{rem:scom}

The bound on the number of iterations in Theorem \ref{thm:GMRES1} crucially depends on the coercivity result of Theorem \ref{thm:coer}. Although numerical experiments in  \cite{BeSp:11} indicate that $\opA$ is coercive, uniformly in $k$, for a wider class of obstacles that those in Theorem \ref{thm:coer}, this has yet to be proved.\footnote{We note that \cite[Remark 6.6]{ChSpGiSm:17} gives an example of a nontrapping obstacle for which 
$\opA$ is \emph{not} coercive uniformly in $k$; therefore, the class of obstacles for which $\opA$ is coercive, uniformly in $k$, is a proper subset of the class of nontrapping obstacles.}

Nevertheless, there \emph{does} exist an integral operator that (i) can be used to solve the sound-soft scattering problem, and (ii) is provable coercive for a wide class of obstacles. Indeed, the \emph{star-combined operator} $\starA$,  introduced in \cite{SpChGrSm:11} and defined by 
\beqs
\starA:= \big(\bx\cdot \bn(\bx)\big)\left(\half I + D'_k\right) + \bx\cdot \nabla_\bound S - \ri \eta S_k
\eeqs
(where $\nT$ is the surface gradient operator on $\bound$; see, e.g., \cite[Page 276]{ChGrLaSp:12}), has the following two properties:
(i) if $u$ solves the sound-soft scattering problem, then 
\beq\label{eq:star_eqn}
\starA \dnpu = \bx\cdot \gamma^+(\nabla u^I)- \ri \eta \gamma^+ u^I
\eeq
\cite[Lemma 4.1]{SpChGrSm:11} (see also \cite[Theorem 2.36]{ChGrLaSp:12}), and 

(ii) if $\Oi$ is a 2- or 3-d Lipschitz obstacle that is star-shaped with respect to a ball and  $\eta := k|\bx| + \ri (d-1)/2$, then
\beqs
\Re\big(\starA \phi,\phi)_{\LtG} \geq \half \essinf_{\bx\in \bound}\big(\bx\cdot\bn(\bx)\big)>0
\eeqs
for all $k>0$ \cite[Theorem 1.1]{SpChGrSm:11}.

The refinement of the Elman estimate in Theorem \ref{thm:Elman} can therefore be used to prove results about the number of iterations required when GMRES is applied to the Galerkin discretisation of \eqref{eq:star_eqn}. 
Since the coercivity constant of the star-combined operator is independent of $k$, the $k$-dependence of the analogue of the bound \eqref{eq:boundonm} for $\starA$ rests on the bounds on $\|\starA\|_{\LtGt}$. 

For convex $\Oi$ with smooth and curved $\bound$, Theorem \ref{thm:L2H1} implies that $\|\starA\|_{\LtGt}\lesssim k^{1/3}$, and we therefore obtain the same bound on $m$ as for $\opA$ (i.e.~\eqref{eq:boundonm}). For general piecewise-smooth Lipschitz obstacles that are star-shaped with respect to a ball, Theorem \ref{thm:L2H1} combined with the bounds \eqref{eq:SbL2} and \eqref{eq:Sb} shows that  $\|\starA\|_{\LtGt}\lesssim k^{1/2}$ when $d=2$ and $\lesssim k^{1/2}\log k$ when $d=3$. Corollary \ref{cor:Elman} then implies that
 $m\gtrsim k^{1/2}$ for $d=2$ and  $m\gtrsim k^{1/2}\log k$ for $d=3$.
Recall that GMRES always converges in at most $N$ steps (in exact arithmetic), and when $h\sim 1/k$ we have that $N\sim k^{d-1}$; these bounds on $m$ are therefore nontrivial.
\ere

\section{Numerical experiments concerning Q2}\label{sec:num}

The main purpose of this section is to show that the $k^{1/3}$ growth in the number of iterations given by Theorem \ref{thm:GMRES1} is effectively sharp.

\paragraph{Details of the scattering problems considered}

We solve the sound-soft scattering problem of Definition \ref{def:SSSP} with $\ba= (1,0,0)$ (i.e the incident plane wave propagates in the $x_1$-direction), using the direct integral equation \eqref{eq:direct} and the Galerkin method \eqref{eq:Galerkin}.
The subspace $\cV_h$ is taken to be piecewise constants on a shape regular mesh, and the meshwidth $h$ is taken to be $2\pi/(10 k)$, i.e.~we are choosing ten points per wavelength. 
We solve the resulting linear system with GMRES, with tolerance $1\times 10^{-5}$.
We consider two obstacles:
\ben
\item $\Oi$ the unit sphere, and
\item $\Oi$ the ellipsoid with semi-principal axes of lengths $3$, $1$, and $1$ (in the $x_1$-, $x_2$-, and $x_3$-directions respectively.
\een
The computations were carried out using version 3.0.3 of the BEM++ library \cite{SmBeArPhSc:15}
on one node of the ``Balena" cluster at the University of Bath. The
cluster consists of Intel Xeon E5-2650 v2 (Ivybridge, 2.60 GHz) CPUs and
the used node had 512GB of main memory. BEM++ was compiled with version
5.2 of the GNU C compiler and the Python code was run under Anaconda
2.3.0.

\paragraph{Numerical results}

\begin{figure}[h!]
  \centering
  \scalebox{.50}
  {
  \includegraphics[width=\textwidth]{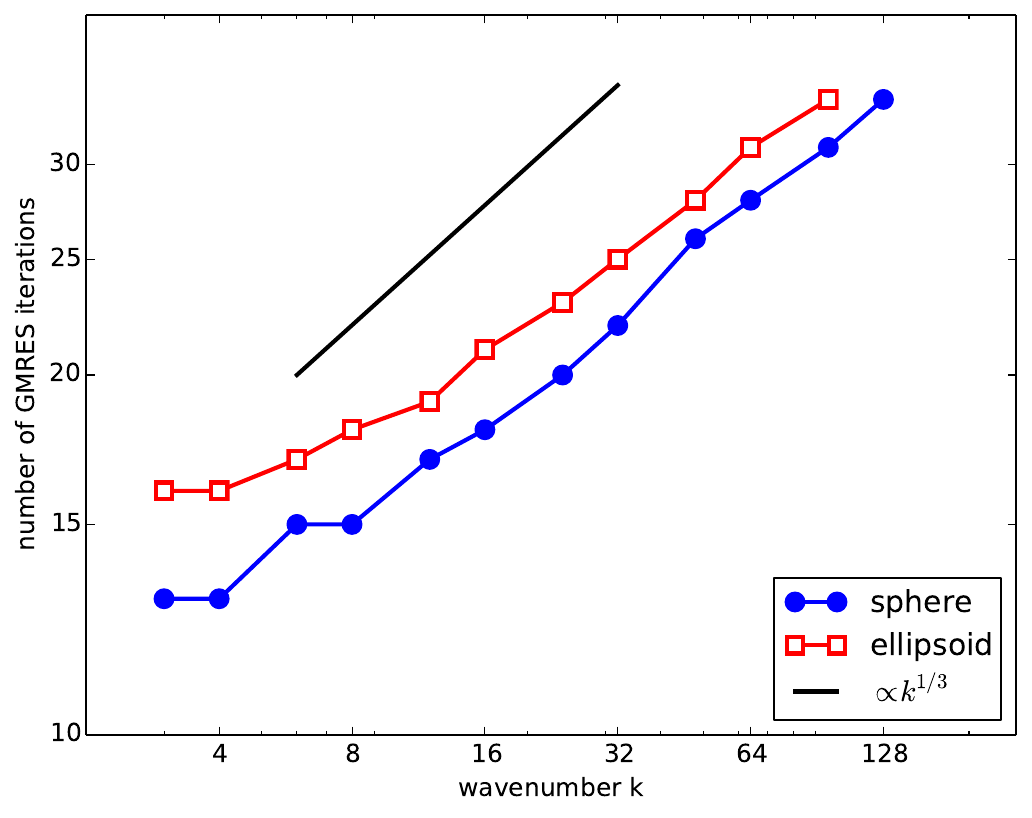}
  }
\caption{The number of iterations required for GMRES to converge (with tolerance $1\times 10^{-5}$) when GMRES is applied to the Galerkin matrix corresponding to the direct integral equation \eqref{eq:direct} with $\eta=k$, and with $\Oi$ the sphere or ellipsoid, and the values of $k$ from Tables \ref{tab:1} and \ref{tab:2}. The $k^{1/3}$ rate is the upper bound on the rate guaranteed by Theorem \ref{thm:GMRES1}.}
\label{fig:1}
\end{figure}

Tables \ref{tab:1} and \ref{tab:2} displays the number of degrees of freedom, number of iterations required for GMRES to converge, and time taken to converge, with $\eta=k$, and with $\Oi$ the sphere or ellipsoid. The difference between Tables \ref{tab:1} and \ref{tab:2} is that, in the first, $k$ starts as $2$ and then doubles until it equals $128$, and in the second, 
$k$ starts as $3$ and then doubles until it equals $96$; we performed the second set of experiments when the $k=128$ run for the ellipsoid failed to complete.
 Figure \ref{fig:1} plots the iteration counts from both tables and compares them to the $k^{1/3}$ rate from Theorem \ref{thm:GMRES1} (the graph is plotted on a $\log$-$\log$ scale so that a dependence $\#_{\text{iterations}}\sim k^\alpha$ appears as a straight line with gradient $\alpha$).
 
\begin{table}[h]
\begin{center}

\begin{tabular}{|c||c|c|c|c|}
\hline 
& \multicolumn{3}{|c|}{Sphere}\tabularnewline
\hline 
$k$ & $\#_{\text{DOF}}$ & $\#_{\text{iterations}}$ & time (s) \tabularnewline
\hline 
 4           &    1304     &    13    &       3.10			\tabularnewline
8         &      4998    &     15   &        7.42	\tabularnewline
16        &      19560    &     18      &    40.30	\tabularnewline
32       &      77224     &    22      &   271.42	\tabularnewline
64       &      307454    &     28     &   2674.54	\tabularnewline
128      &      1225260   &      34    &   31024.43	\tabularnewline
\hline 
\end{tabular}~%
\begin{tabular}{|c|c|c|c|}
\hline 
\multicolumn{3}{|c|}{Ellipsoid}\tabularnewline
\hline 
$\#_{\text{DOF}}$ & $\#_{\text{iterations}}$ & time (s) \tabularnewline
\hline 
3230      &   16       &    5.26 			\tabularnewline
   12324  &       18     &     19.30		\tabularnewline
   48526    &     21     &    113.95		\tabularnewline
       190784  &       25    &     926.47	\tabularnewline
    754236     &    31  &     10354.29	\tabularnewline
   \hspace{0.01ex} * & *& *\tabularnewline
\hline 
\end{tabular}
\end{center}
\caption{
With $\Oi$ the sphere or ellipsoid and $\eta=k$,
the number of degrees of freedom, number of iterations required for GMRES to converge (with tolerance $1\times 10^{-5}$), and time taken to converge, when GMRES is applied to the Galerkin matrix corresponding to the direct integral equation \eqref{eq:direct}, starting with $k=4$ and then doubling until $k=128$. $^*$ denotes that the run did not complete.
}\label{tab:1} 
\end{table}

\begin{table}[h]
\begin{center}

\begin{tabular}{|c||c|c|c|c|}
\hline 
& \multicolumn{3}{|c|}{Sphere}\tabularnewline
\hline 
$k$ & $\#_{\text{DOF}}$ & $\#_{\text{iterations}}$ & time (s) \tabularnewline
\hline 
3             &   846       &  13        &   1.12	\tabularnewline
6              & 2880       &  15         &  3.85	\tabularnewline
12             & 11054      &   17        &  18.56	\tabularnewline
24             & 43688       &  20       &  107.18	\tabularnewline
48            & 173264       &  26       &  928.61	\tabularnewline
96            & 689894       &  31      & 10753.95	\tabularnewline
\hline 
\end{tabular}~%
\begin{tabular}{|c|c|c|c|}
\hline 
\multicolumn{3}{|c|}{Ellipsoid}\tabularnewline
\hline 
$\#_{\text{DOF}}$ & $\#_{\text{iterations}}$ & time (s) \tabularnewline
\hline 
   1806     &    16       &    6.20		\tabularnewline
   6874       &  17       &    9.51		\tabularnewline
  26994      &   19      &    55.64		\tabularnewline
 107272      &   23     &    373.45	\tabularnewline	
  426026      &   28     &   3985.63	\tabularnewline
  1691328      &   34     &  43423.69	\tabularnewline
\hline 
\end{tabular}
\end{center}

\caption{ Same as Table \ref{tab:1} but for a different range of $k$.
}\label{tab:2} 
\end{table}

We see from Figure \ref{fig:1} that the $k^{1/3}$ growth predicted by Theorem \ref{thm:GMRES1} appears to be effectively sharp. Indeed, the plot of the iterations for the ellipsoid becomes roughly linear from $k=12$ onwards, and estimating the slope of this line using the numbers of iterations at $k=12$ and $k=96$ we have that the  $\#_{\text{iterations}}\sim k^{0.28}$. 
Using the numbers of iterations at $k=12$ and $k=96$ to estimate the rate of growth for the sphere we have that $\#_{\text{iterations}}\sim k^{0.29}$.

Finally, Table \ref{tab:3} compares the iteration counts and times for the sphere when $\eta=k$ and when $\eta=-k$. We see that, for every value of $k$ considered, the number of iterations when $\eta=-k$ is much greater than when $\eta=k$. Table \ref{tab:3} only goes up to $k=32$, since the $k=64$ run for the sphere with $\eta=-k$ did not complete.

\begin{table}[h]
\begin{center}

\begin{tabular}{|c||c|c|c|c|}
\hline 
& \multicolumn{2}{|c|}{$\eta=k$}\tabularnewline
\hline 
$k$ & $\#_{\text{iterations}}$ & time (s) \tabularnewline
\hline 
4       &              13       &    3.10 	\tabularnewline 	
8          &             15      &     7.42	\tabularnewline
16          &         18       &   40.30	\tabularnewline
32           &          22     &    271.42	\tabularnewline
\hline 
\end{tabular}~%
\begin{tabular}{|c|c|c|c|}
\hline 
\multicolumn{2}{|c|}{$\eta=-k$}\tabularnewline
\hline 
 $\#_{\text{iterations}}$ & time (s) \tabularnewline
\hline 
        44     &      3.46		\tabularnewline
      88        &   9.04		\tabularnewline
     405        &  75.38		\tabularnewline
     11191      &  4502.05	\tabularnewline
\hline 
\end{tabular}
\end{center}
\caption{With $\Oi$ the sphere and $\eta=k$ or $\eta=-k$, the number of iterations required for GMRES to converge (with tolerance $1\times 10^{-5}$) and time taken to converge, when GMRES is applied to the Galerin matrix corresponding to the direct integral equation \eqref{eq:direct}. 
}\label{tab:3} 
\end{table}

\bre[The link between Table \ref{tab:3} and the recent work of Marburg \cite{Ma:14, Ma:15}]

We performed the experiment in Table \ref{tab:3} because, in the engineering-acoustics literature, Marburg recently considered collocation discretisations of the direct integral equation for the Neumann problem (i.e.~the Neumann-analogue of equation \eqref{eq:direct}) and showed that the analogue of the choice $\eta=k$ leads to much slower growth than the analogue of the choice $\eta=-k$ \cite{Ma:14, Ma:15}. 

A heuristic explanation for this dependence of the number of iterations on the sign of $\eta$ is essentially contained in the work of Levadoux and Michielsen
\cite{Le:01, LeMi:04}, and Antoine and Darbas \cite{AnDa:05}. In our setting of using the operator $\opA$ to solve the exterior Dirichlet problem, the key points are that 
\ben
\item the ideal $\ri \eta$ should approximate the Dirichlet-to-Neumann (DtN) map in $\Oe$, and
\item $\ri k$ is a better approximation to the DtN map than $-\ri k$ (at least for smooth convex obstacles).
\een

Regarding 1: taking the Dirichlet trace of Green's integral representation (written with general Dirichlet data, not just data coming from a plane wave as in \eqref{eq:IR}), 
and using the jump relations for the single- and double-layer potentials (see, e.g., \cite[Equation 2.41]{ChGrLaSp:12}) we find that
\beqs
\gpu = -S_k(\dnpu)+ \left(\half I+ D_k\right)\gpu.
\eeqs
Rearranging this equation, and introducing the notation $\DtN^+$ for the exterior Dirichlet-to-Neumann map
for solutions of the Helmholtz equation in $\Oe$ satisfying the Sommerfeld radiation condition, we find that
\beq\label{eq:A1}
I = \half I + D_k- S_k \DtN^+.
\eeq
Green's second identity implies that, for $\phi,\psi \in H^{1/2}(\bound)$, 
\beqs
\langle \DtN^+\phi, \psi \rangle_\bound= \langle \phi, \DtN^+\psi \rangle_\bound,
\eeqs 
where the duality pairing $\langle\phi, \psi\rangle_\bound := \int_\bound \phi\, \psi \,\rd s$ when $\phi, \psi\in \LtG$; see \cite[Equation 2.65, Equation 2.84, Equation A.24]{ChGrLaSp:12}.
Therefore, taking the adjoint of \eqref{eq:A1}, we find that
\beq\label{eq:A1a}
I = \half I + D'_k- \DtN^+ S_k.
\eeq
Comparing \eqref{eq:A1a} to the definition of $\opA$ in \eqref{eq:scpo}, we see that the ideal $\ri \eta$ should approximate $\DtN^+$. 
The idea of choosing $\eta$ as an operator, based on the relations \eqref{eq:A1}-\eqref{eq:A1a} (and their analogues for the Neumann problem), essentially first appeared in \cite{Le:01}, \cite{LeMi:04}. The relation \eqref{eq:A1} appeared explicitly in \cite[Theorem 2.1]{AnDa:05}, with this paper considering local approximations of the non-local operators $\DtN^+$ and $\NtD^+$, whilst \cite{Le:01}, \cite{LeMi:04} used non-local pseudodifferential-operator approximations. 

Regarding 2: In the case when $\partial \Omega$ is a circle (of radius $1$), the DtN map is given by 
\beqs
\pdiff{u}{r}(1,\theta) = k\sum_{n=-\infty}^\infty  \frac{\Hnpka}{\Hna}  \re^{\ri n \theta} \, \left(\frac{1}{2\pi} \intotp \re^{-\ri n \phi}\, u(1,\phi) \, \rd\phi\right).
\eeqs
The uniform- and double-asymptotic expansions of the Hankel functions (see, e.g., \cite[\S10.20, 10.41(v)]{Di:16}) imply that
\beq\label{eq:Hankelkey}
k  \frac{\Hnpka}{\Hna} \sim
\begin{cases}
\ri k,  &\tfor n \text{ fixed as } k \tendi,\\
\ri k \sqrt{1- \left(\frac{n}{k}\right)^2}, & \,\, n, k \tendi \text{ with } k-|n| \gg k^{1/3},\\
\re^{2\pi \ri /3}\sqrt{\frac{n^2-k^2}{n^{2/3}\zeta}} \frac{ {\rm Ai}' \big(\re^{2\pi \ri/3}n^{2/3}\zeta\big)}{{\rm Ai}\big(\re^{2\pi \ri/3}n^{2/3}\zeta\big)},& \,\, n, k \tendi \text{ with } \big||n|-k\big|\leq  Mk^{1/3},\\
n\sqrt{1-\left(\frac{k}{n}\right)^2}, & \,\, n, k \tendi \text{ with } |n|-k\gg k^{1/3},
\end{cases}
\eeq
where $\zeta$ is defined in terms of $n$ and $k$ by \cite[Equations 10.20.2 and 10.20.3]{Di:16} \footnote{Strictly speaking, the case $n/k\tendi$ is not covered in the asymptotics \cite[\S10.20, 10.41(v)]{Di:16}, but the asymptotics here can be shown using more general microlocal methods \cite{Ga:15a}.}.
We see that the approximation $k  \Hnpka/\Hna\sim \ri k$
describes the DtN map on the low frequency modes and, in particular, is much better than the approximation $k  \Hnpka/\Hna\sim -\ri k$. 
The asymptotics \eqref{eq:Hankelkey}, however, show that neither the approximations $\ri k$ or $-\ri k$ are particularly good on the higher frequency modes. 
An almost-identical analysis is valid for the sphere, and more generally for a smooth convex curved obstacle, since the symbol of the DtN map for such domains is described by the asymptotics \eqref{eq:Hankelkey}; see \cite[\S9, last formula on page 58]{Ga:15a}.
\ere
 
 \paragraph{Acknowledgements.}
The authors thank Timo Betcke (University College London), Simon Chandler-Wilde (University of Reading), and S\'ebastien Loisel (Heriot-Watt University) for useful discussions.
In particular, EAS's understanding of the improvement on the Elman estimate in \cite{BeGoTy:06} was greatly improved after seeing the alternative derivation of this result in \cite[Equations 3.2 and 4.8]{GrLo:15}.
This research made use of the Balena High Performance Computing (HPC) Service at the University of Bath.  JG thanks the US National Science Foundation for support under the Mathematical Sciences Postdoctoral Research Fellowship  DMS-1502661, and EAS thanks the UK Engineering and Physical Sciences Research Council for support under Grant EP/R005591/1.

\def\cprime{$'$} \def\cprime{$'$} \def\cprime{$'$}
\footnotesize{
\bibliographystyle{plain}
\bibliography{biblio_acta}
}

\end{document}